\documentstyle[amssymb, 12pt]{article}

\begin{document}
\font\frak=eufm10 scaled\magstep1
\font\fak=eufm10 scaled\magstep2
\font\fk=eufm10 scaled\magstep3
\font\scriptfrak=eufm10
\font\tenfrak=eufm10

\newtheorem{theorem}{Theorem}
\newtheorem{corollary}{Corollary}
\newtheorem{proposition}{Proposition}
\newtheorem{definition}{Definition}
\newtheorem{lemma}{Lemma}
\font\frak=eufm10 scaled\magstep1
\newenvironment{pf}{{\noindent{\it Proof. }}}{\ $\Box$\medskip}


\mathchardef\za="710B  
\mathchardef\zb="710C  
\mathchardef\zg="710D  
\mathchardef\zd="710E  
\mathchardef\zve="710F 
\mathchardef\zz="7110  
\mathchardef\zh="7111  
\mathchardef\zvy="7112 
\mathchardef\zi="7113  
\mathchardef\zk="7114  
\mathchardef\zl="7115  
\mathchardef\zm="7116  
\mathchardef\zn="7117  
\mathchardef\zx="7118  
\mathchardef\zp="7119  
\mathchardef\zr="711A  
\mathchardef\zs="711B  
\mathchardef\zt="711C  
\mathchardef\zu="711D  
\mathchardef\zvf="711E 
\mathchardef\zq="711F  
\mathchardef\zc="7120  
\mathchardef\zw="7121  
\mathchardef\ze="7122  
\mathchardef\zy="7123  
\mathchardef\zf="7124  
\mathchardef\zvr="7125 
\mathchardef\zvs="7126 
\mathchardef\zf="7127  
\mathchardef\zG="7000  
\mathchardef\zD="7001  
\mathchardef\zY="7002  
\mathchardef\zL="7003  
\mathchardef\zX="7004  
\mathchardef\zP="7005  
\mathchardef\zS="7006  
\mathchardef\zU="7007  
\mathchardef\zF="7008  
\mathchardef\zW="700A  

\newcommand{\be}{\begin{equation}}
\newcommand{\ee}{\end{equation}}
\newcommand{\ra}{\rightarrow}
\newcommand{\lra}{\longrightarrow}
\newcommand{\bea}{\begin{eqnarray}}
\newcommand{\eea}{\end{eqnarray}}
\newcommand{\beas}{\begin{eqnarray*}}
\newcommand{\eeas}{\end{eqnarray*}}
\newcommand{\Z}{{\Bbb Z}}
\newcommand{\R}{{\Bbb R}}
\newcommand{\C}{{\Bbb C}}
\newcommand{\1}{{\bold 1}}
\newcommand{\SL}{SL(2,\C)}
\newcommand{\Sl}{sl(2,\C)}
\newcommand{\SU}{SU(2)}
\newcommand{\su}{su(2)}
\newcommand{\G}{{\goth g}}
\newcommand{\D}{{\rm d}}
\newcommand{\de}{\,{\stackrel{\rm def}{=}}\,}
\newcommand{\we}{\wedge}
\newcommand{\nn}{\nonumber}
\newcommand{\ot}{\otimes}
\newcommand{\s}{{\textstyle *}}
\newcommand{\ts}{T^\s}
\newcommand{\da}{\dagger}
\newcommand{\pa}{\partial}
\newcommand{\ti}{\times}
\newcommand{\A}{{\cal A}}
\newcommand{\Li}{{\cal L}}
\newcommand{\ka}{{\Bbb K}}
\newcommand{\find}{\mid}
\newcommand{\ad}{{\rm ad}}
\newcommand{\rS}{]^{SN}}
\newcommand{\rb}{\}_P}
\newcommand{\p}{{\sf P}}
\newcommand{\h}{{\sf H}}
\newcommand{\X}{{\cal X}}
\newcommand{\I}{\,{\rm i}\,}
\newcommand{\rB}{]_P}
\newcommand{\Ll}{{\pounds}}

\title{Jacobi structures revisited}

\author{
Janusz Grabowski\thanks{Supported by KBN, grant No. 2 P03A 031 17.}\\
Institute of Mathematics, Warsaw University\\
ul. Banacha 2, 02-097 Warszawa, Poland. \\
and\\
Mathematical Institute, Polish Academy of Sciences\\
ul. \'Sniadeckich 8, P. O. Box 137, 00-950 Warszawa, Poland\\
{\it e-mail:} jagrab@mimuw.edu.pl
\and
Giuseppe Marmo\thanks{Supported by PRIN SINTESI}\\
Dipartimento di Scienze Fisiche,
Universit\`a Federico II di Napoli\\
and\\
INFN, Sezione di Napoli\\
Complesso Universitario di Monte Sant'Angelo\\
Via Cintia, 80126 Napoli, Italy\\
{\it e-mail:} marmo@na.infn.it}
\maketitle
\begin{abstract} Jacobi algebroids, that is graded Lie brackets
on the Grassmann algebra associated with a vector bundle which satisfy a
property similar to that of Jacobi brackets are  introduced.  They  turn
out to be equivalent to generalized  Lie  algebroids  in  the  sense  of
Iglesias and Marrero. Jacobi  bialgebroids  are  defined  in  the  same
manner.
A lifting procedure of elements of this Grassmann algebra to multivector
fields on the total space of the  vector  bundle  which  preserves  the
corresponding Lie brackets is developed. This gives the possibility  of
associating canonically a Lie algebroid with any local Lie algebra in
the sense of Kirillov.
\end{abstract}

\section{Introduction}
This work was originated as an attempt to understand  the  Lie  algebroid
structure  on  $T^*M\oplus_M\R$  associated  with  a  Jacobi   structure
$(\zL,\zG)$ on a manifold $M$ (\cite{Li2}). The formula in \cite{KSB}
\bea\nn
[(\za,f),(\zb,g)]&=&(\Li_{\zL_\za}\zb-\Li_{\zL_\zb}\za-\D<\zL,\za\we\zb>
+f\Li_\zG\zb-g\Li_\zG\za-i_\zG\za\we\zb,\label{jb}\\ \nn
&&<\zL,\zb\we\za>+\zL_\za(g)-\zL_\zb(f)+f\zG(g)-g\zG(f)),
\eea
which gives the Lie bracket on the space  $\zW^1(M)\ti  C^\infty(M)$  of
sections of  $T^*M\oplus_M\R$,  being  rather  complicated,  deserves  a
better understanding and explanation. During our work  we  have  noticed
that it is very close to \cite{IM} which will be our  primary  reference
paper.

Since a Jacobi bracket is just a Lie bracket on the  algebra  of  smooth
functions given by a bilinear first-order differential operator, we  start
with the study of the Nijenhuis-Richardson bracket on multilinear first order
differential operators. This bracket is a  graded  Lie  bracket  but  it
differs from the Richardson-Nijenhuis bracket. This is manifested  by  the
fact that with respect to the wedge product it is not a derivation
but  a first order differential operator. This is like  the  difference
between Poisson and Jacobi brackets.

We then discuss the case of a general Lie
algebroid.  Our  primary  object  is  the   Schouten-Nijenhuis   bracket
associated with the Lie algebroid rather than the Lie  algebroid  bracket
itself. Deforming the Schouten-Nijenhuis bracket to a graded Lie bracket
which violates the   Leibniz   rule, like   in   the   case   of   the
Nijenhuis-Richardson bracket for first order differential operators, we
introduce the notion of a {\it Jacobi algebroid}. We find out that  this
is a structure equivalent to the notion of a {\it  Lie  algebroid   with
the presence of 1-cocycle} as defined in \cite{IM}.

Since any Lie algebroid structure on a vector bundle $E$  is  associated
with a linear Poisson structure on the dual bundle $E^*$, one can expect
that  there  is  a  lifting  procedure   of   multilinear   first-order
differential operators acting on smooth functions on $M$ to  multivector
fields on $TM\oplus_M\R$, similar to the classical complete tangent lift
of multivector fields on $M$ (cf. \cite{IY,GU}), which associates 
the corresponding linear Poisson structure with a given Jacobi bracket.
We define such lifts for Jacobi algebroids and show that the lift of a
Jacobi structure gives exactly the Lie algebroid bracket (\ref{br}).
We extend this for general local Lie algebra structure in  the
sense of Kirillov \cite{Ki}. The main  result  is  that  any  local  Lie
algebra structure on a one-dimensional bundle $L$ induces naturally  a
Lie algebroid structure on the first jet bundle $J_1(L)$.

Introducing a Cartan calculus for a given Jacobi  Lie  algebroid as in
\cite{IM} one can define {\it Jacobi bialgebroids}, by  analogy  to  Lie
bialgebroids, as Jacobi algebroid  structures  on  dual  pair  of  vector
bundles such that the exterior differential induced by one structure  is
a graded derivation for the Schouten-Jacobi bracket of the  second  one.
We show that this reduces  exactly to the notion  of  a  {generalized  Lie
bialgebroid} in \cite{IM}. The advantage of  using  consequently  graded
brackets on the corresponding Grassmann algebras is that this  definition
becomes more natural.

\setcounter{equation}{0}
\section{Graded Lie brackets}

A  {\em  graded  Lie  bracket\/}  on   a   graded   vector   space
$\A=\bigoplus_{n\in\Z}\A^n$  (`graded' means   always  `$\Z$-graded'
throughout this paper) is a bilinear operation $[\  ,\  ]:\A\times
\A\lra \A,$ being graded
\be
[\A^n,\A^m]\subset \A^{n+m},
\ee
graded skew-symmetric
\be
[X,Y]=-(-1)^{xy}[Y,X],
\ee
and satisfying the graded Jacobi identity
\be
[[X,Y],Z]=[X,[Y,Z]]-(-1)^{xy}[Y,[X,Z]],
\label{Jacobi}
\ee
where we fix the convention that we write simply $x,a$,  etc.,  for  the
Lie algebra degrees of homogeneous elements $X,A$, etc., when  no
confusion arises.

One sometimes writes the graded Jacobi identity in the form
\be
(-1)^{xz}[[X,Y],Z]+(-1)^{yx}[[Y,Z],X]+(-1)^{zy}[[Z,X],Y]=0
\label{Jacobi1}
\ee
which is equivalent to (\ref{Jacobi}) for  graded  skew-symmetric
brackets. However, for  non-skew-symmetric  brackets  the  formula
(\ref{Jacobi}) seems to  be  better,  since  it  means  that  the
adjoint map $X\mapsto \ad_X\de [X,\cdot ]$ is a representation  of
the bracket, i.e. $\ad_{[X,Y]}$ is equal to the graded commutator
\be
[\ad_X,\ad_Y]\de\ad_X\circ\ad_Y-(-1)^{xy}\ad_Y\circ\ad_X=
\ad_{[X,Y]},
\label{ad}
\ee
whereas (\ref{Jacobi1}) has no clear direct meaning.
\par
With a given smooth ($C^\infty$) manifold $M$ 
several natural graded Lie brackets of tensor fields are associated.
Historically
the first  one  was  probably  the  famous  Schouten-Nijenhuis   bracket
$[\cdot,\cdot\rS$ defined on  multivector  fields  (see  \cite{Sc,Ni}).
It is the unique graded extension of the usual bracket $[\cdot,\cdot]$
on the space $\X (M)$ of vector  fields  to  the  exterior  algebra
$A(M)=\bigoplus_{n\in\Z}A^{[n]}(M)$  of  multivector   fields (where
$A^{[n]}(M)=\zG (\zL^nTM)$ is the space of $n$-vector fields for $n\geq
0$ and $A^{[n]}(M)=\{ 0\}$ for $n<0$) such that
\begin{description}
\item{(a)} the degree of $X\in A^{[n]}(M)$ with respect to the   bracket
is $(n-1),$
\item{(b)}
$[X,f\rS  =X(f)$,
\item{(c)}
$ [X,Y\wedge Z\rS =[X,Y\rS\wedge Z+(-1)^{(k-1)l}Y\wedge[X,Z\rS,\qquad$
for $X\in  A^{[k]}(M),\,  Y\in  A^{[l]}(M),$  i.e.  $\ad$  is  a
representation of the Schouten-\-Nijenhuis bracket in graded
derivations  of  the  graded associative algebra $A(M).$
\end{description}
The Schouten-Nijenhuis bracket is an example of what is sometimes called
a {\em Gerstenhaber  algebra\/}   (see   \cite{KS,KS1}) which
consists  of  a  triple
$(\A=\oplus_{n\in\Z}\A^n,\we,[\cdot,\cdot])$,          such          that
$(\A,[\cdot,\cdot])$     is     a     graded     Lie     algebra     and
$(\A=\oplus_{n\in\Z}\A^{[n]},\we)$, with $\A^{[n]}=\A^{n+1}$, is a graded
associative commutative algebra, and $\ad_X$ for $X\in\A^x$ is a
derivation with
respect to $\we$ with degree $x$, i.e.
\be    [X,Y\wedge    Z]    =[X,Y]\wedge     Z+(-1)^{x(y+1)}Y\wedge[X,Z].
\label{Schouten2}
\ee
From (\ref{Schouten2}) it follows that in any Gerstenhaber algebra
\begin{eqnarray}\label{Schouten3}
\lefteqn
{[X_1\wedge\ldots\wedge X_m,Y_1\wedge\cdots\wedge Y_n] =} \\
& & \sum_{k,l}(-1)^{k+l}[X_k,Y_l]\wedge\ldots\wedge\widehat{X_k}\wedge
\ldots\we X_m\we Y_1\we\ldots\we\widehat{Y_l}\we\ldots\we Y_n,
\nonumber
\end{eqnarray}
where $X_k,Y_l\in\A^0$ and the  hats  stand  for  omissions.  Note  that
$\A^0$ is a Lie subalgebra of $\A$ and  $V=\A^{-1}$  is  an  associative
commutative subalgebra of $\A$.
\par
The Schouten-Nijenhuis bracket is a particular case  of  the  so-called
{\it   Nijenhuis-Richardson   bracket} (\cite{NR})  which,   in    turn,
is    the skew-symmetrization of the natural graded Lie bracket on
multilinear operators  discovered by Gerstenhaber \cite{Ge}. We  will
call  the last  bracket  the  {\it Gerstenhaber bracket}. We will define
all these brackets in  details  to fix notation and signs.

\medskip\noindent
{\it Gerstenhaber bracket.} Let $V$ be  a  vector  space  over  a  field
$\frak  k$.  Denote  by  $M^p(V)$  the  space  of  $(p+1)$-linear   maps
$A:V\ti\dots\ti    V\ra    V$.    On    the    graded    vector    space
$M(V)=\bigoplus_{p\in\Z}M^p(V)$, where $M^{-1}(V)=V$ and $M^p(V)=\{ 0\}$
for $p<-1$, we define first {\it insertion operators}. For $A\in M^a(V),
B\in M^b(V)$, the insertion $\i_BA\in M^{a+b}(V)$ is defined by
\be
\i_BA(x_0,\dots,x_{a+b})=\sum_{k=0}^a(-1)^{bk}A(x_0,\dots,x_{k-1},
B(x_k,\dots,x_{k+b}),x_{k+b+1},\dots,x_{a+b}).
\ee
Then the {\it Gerstenhaber bracket} $[A,B]^G$ is given by
\be
[A,B]^G=-\i_AB+(-1)^{ab}\i_BA.
\ee
We  can  consider  the  graded  subspace   $\A(V)$   of   skew-symmetric
elements of $M(V)$. The {\it Nijenhuis-Richardson bracket} on $\A(V)$  is
the skew-symmetrization of the Gerstenhaber bracket:
\be
[A,B]^{NR}=\frac{(a+b+1)!}{(a+1)!(b+1)!}{\rm skew}([A,B]^G),
\ee
where ${\rm skew}$ stands for the  antisymmetrization  projector
in $M(V)$. We have $[A,B]^{NR}=-i_AB+(-1)^{ab}i_BA$, with
\be
i_BA(x_0,\dots,x_{a+b})=\sum_{\zs\in  S(a+b,b)}(-1)^\zs
A(B(x_{\zs(0)},\dots,x_{\zs(b)}),x_{\zs(b+1)},\dots,x_{\zs(a+b)}),
\ee
where $S(a+b,b)$ is the set of unshuffles $\zs:\{ 0,\dots,a+b\}\ra
\{ 0,\dots,a+b\}$ with $\zs(0)<\dots<\zs(b)$, $\zs(b+1)<\dots<\zs(a+b)$.
The following is well known.
\begin{theorem} The Gerstenhaber  bracket  and  the  Nijenhuis-Richardson
bracket make the graded vector spaces $M(V)$ and $\A(V)$,  respectively,
into graded Lie algebras. Moreover,
\item{a)} the equation $[A,A]^G=0$ for $A\in M^1(V)$  is  equivalent  to
the fact that the bilinear operation $A$ on $V$ is associative;
\item{b)} the equation $[A,A]^{NR}=0$ for $A\in \A^1(V)$  is  equivalent
to the fact that the  skew-bilinear  operation  $A$  on  $V$  is  a  Lie
bracket.
\medskip
\end{theorem}
Suppose now that $V$ is an associative  commutative  algebra  with  unit
$\1$. Then $M(V)$ is a graded associative algebra (with elements  of
$M^a(V)$ being of degree $a+1$) with the obvious product
\be
A\cdot
B(x_0,\dots,x_{a+b+2})=A(x_0,\dots,x_a)B(x_{a+1},\dots,x_{a+b+2}),
\ee
for $A\in\A^a(V)$, $B\in\A^b(V)$. Similarly, $\A(V)$ is in a natural way
a graded associative commutative algebra (again, with elements  of
$M^a(V)$ being of degree $a+1$) with the obvious wedge product
\be
A\we B=\frac{(a+b+2)!}{(a+1)!(b+1)!}{\rm skew}(A\cdot B).
\ee
Denote  by   $Diff(V)$,   $Diff_1(V)$,   and   $Der(V)$,   respectively,
the  space  of  linear  differential   operators,   linear   first order
differential operators,  and  derivations  on  $V$.  Similarly,  by  $\A
Diff^p(V)$,
$\A  Diff^p_1(V)$,  and  $\A  Der^p(V)$,  we  denote  the  corresponding
skew  $(p+1)$-linear  operators  on  $V$   which   are,    respectively,
differential  operators  ,  first order  differential   operators,   and
derivations with respect to each  variable separately. By $\A Diff(V)$,
$\A Diff_1(V)$, and $\A Der(V)$, we  denote the corresponding graded
vector spaces. It is easy to see  the  following (cf. \cite{Gr}):
\begin{theorem}
\item{(a)} $\A Diff(V)$, $\A Diff_1(V)$,  and  $\A  Der(V)$  are   graded
associative subalgebras of $(\A(V),\we)$ and graded  Lie subalgebras of
$(\A(V),[\cdot,\cdot]^{NR})$;
\item{(b)} There is a canonical splitting
\be
\A Diff^p_1(V)=\A Der^p(V)\oplus \A Der^{p-1}(V)
\ee
given   by   $A=A_1+I\we   A_2$,   where   $A_1\in\A     Der^p(V)$,
$A_2=i_\1A=A(\1,\cdot,\dots,\cdot)\in\A Der^{p-1}(V)$, and $I$ is the
identity map on $V$.
\item{(c)}  $(\A  Der(V),\we,[\cdot,\cdot]^{NR})$  is   a   Gerstenhaber
algebra. In the case when  $V=C^\infty(M)$  is  the  algebra  of  smooth
functions  on  a  manifold  $M$,  one  has  $\A  Der(V)=A(M)$  and   the
Nijenhuis-Richardson bracket reduces to the Schouten-Nijenhuis bracket.
\medskip
\end{theorem}
\begin{pf} One can find the proofs of the parts (a) and (b) 
in \cite{Gr}, Section 3. The part (c) follows easily from the  following
properties of the insertion operators versus wedge products:
\bea\label{we1}
i_C(A\we B)&=&(i_CA)\we B+(-1)^{c(a+1)}A\we(i_CB),\\
i_{A\we B}C&=&(-1)^{c(b+1)}(i_AC)\we B+A\we(i_BC)\quad {\rm for}\quad
C\in\A Der^c(V).\label{we2}
\eea
Here, according to our convention, $A\in  \A^a(V)$,  i.e.  $A$  is
$(a+1)$-linear, etc.
\end{pf}

{\bf Remark.} Our convention of signs in the Gerstenhaber, and hence  in
the Nijenhuis-Richardson bracket, is different  from  the  original  one.
This is chosen in this way  in  order  to  get  a  Gerstenhaber  algebra
structure  on $\A Der(V)$ and  hence  on  $A(M)$.  Also  the  standard
Schouten-Nijenhuis bracket, which is still used by many  authors,
differs  by  sign   from our.
In particular, the Schouten bracket used in \cite{IM}
is not a graded Lie bracket, since it is not graded  skew-symmetric.  It
seems reasonable to use consequently  graded  Lie  algebra  brackets  in
order  to  avoid   confusions.   This   will  also   simplify    certain
formulae and definitions, as we will see it later.

\medskip
Since we already know that for an associative  commutative  algebra  $V$
the triple $(\A Der(V),\we,[\cdot,\cdot]^{NR})$ is a Gerstenhaber
algebra, let us  look closer at the structure of the algebra $(\A
Diff_1(V),\we,[\cdot,\cdot]^{NR})$.  In  the
case $V=C^\infty(M)$  we  will  write  $\A  Diff_1(M)$  instead  of  $\A
Diff_1(V)$.
\par
Since $i_{A\we B}I=A\we B=(i_AI)\we B+A\we(i_BI)-A\we  B$,
we have  for  $C\in\A Diff_1^c(V)$  instead  of  (\ref{we2})  the
following:
\be
i_{A\we B}C=(-1)^{c(b+1)}(i_AC)\we B+A\we(i_BC)-A\we B\we i_\1C.
\ee
This, in turn, implies that on $\A Diff_1(V)$ we have
\be\label{GJ}
[A,B\we C]^{NR}=[A,B]^{NR}\we C+(-1)^{a(b+1)} B\we[A,C]^{NR}
-(-1)^ai_\1A\we B\we C.
\ee
Note that $D=i_\1$ is a graded derivative of the wedge product of degree -1
and $\tilde D(X)=(-1)^xD(X)$ defines a right derivative:
\be
\tilde D(X\we Y)=X\we\tilde D(Y)+(-1)^{y+1}\tilde D(X)\we Y.
\ee
In general, we will call a {\it Gerstenhaber-Jacobi algebra} a triple
$(\A=\oplus_{n\in\Z}\A^n,\we,[\cdot,\cdot])$  as  in  the   Gerstenhaber
algebra case but with the graded bracket satisfying
\be\label{F}
[X,Y\we Z]=[X,Y]\we Z+(-1)^{x(y+1)}Y\we[X,Z]-\tilde D(X)\we Y\we Z,
\ee
where $\tilde D$ is a graded linear map of degree -1, instead  of  the
Leibniz rule.
Putting $Y=Z=\1$ we obtain that $\tilde D(X)=[X,\1]$, so that $\tilde  D$
is a right graded derivation of degree -1  with  respect  to  both:  the
associative and Lie algebra structures.
Here we assume that the associative commutative algebra  $V=\A^{-1}$
has the unit $\1$  (if  not,  we  can  always  easily  extend  the  whole
structure).
Thus, in  the   case of a Gerstenhaber-Jacobi algebra $\ad_X$   is   not
a  derivative but  a differential operator of first order with respect to
the  wedge  product (cf. \cite{Ko}).

\section{Lie algebroids and Jacobi algebroids}
Let $M$ be a smooth manifold. A {\em Lie algebroid  on
$M$} is a vector bundle $\zt:L\ra M$, together with a  bracket
$[\cdot,\cdot]:\zG L\ti \zG L\ra\zG L$ on the $C^\infty(M)$-module  $\zG
L$ of smooth sections of $L$, and a $C^\infty(M)$-linear map $a:\zG  L\ra
\X(M)$ from $\zG L$ to the Lie algebra of vector fields on $M$, called
the {\em  anchor}  of the Lie algebroid, such that
\noindent
\begin{description}
\item{(i)} the bracket on $\zG L$ is $\R$-bilinear, alternating, and
satisfies the Jacobi identity;
\item{(ii)} $[X,fY]=f[X,Y]+a(X)(f)Y$ for all  $X,Y\in  \zG  L$  and  all
$f\in C^\infty(M)$.
\end{description}
From (i) and (ii) it follows easily
\begin{description}
\item{(iii)} $a([X,Y])=[a(X),a(Y)]$ for all $X,Y\in\zG L$.
\end{description}
We get an algebraic counterpart of the notion of Lie algebroid replacing
the  algebra  $C^\infty(M)$  of  smooth  functions   by   an   arbitrary
associative commutative algebra $V$, and the module of sections  of   the
vector bundle $\zt:L\ra M$ by a module $\Li$ over the algebra $V$: a {\em
Lie pseudoalgebra} over  $V$ is a $V$-module $\Li$ together with a
bracket
$[\cdot,\cdot]:\Li\ti\Li\ra\Li$ on the module $\Li$, and a $V$-module
morphism
$a:\Li\ra{\rm Der}(V)$ from $\Li$ to the $V$-module ${\rm Der}(V)$ of
derivations of $V$, called the {\em anchor} of $\Li$, such that
\begin{description}
\item{(i)} the bracket on $\Li$ is bilinear , alternating, and
satisfies the Jacobi identity;
\item{(ii)} For all  $X,Y\in  \Li$  and  all $f\in V$ we have
\be\label{0}
[X,fY]=f[X,Y]+a(X)(f)Y;
\ee
\item{(iii)} $a([X,Y])=[a(X),a(Y)]$ for all $X,Y\in\Li$.
\end{description}
As before, (i) and (ii) imply (iii) if  only  the  $V$-module  $\Li$  is
faithful.
\par
Lie algebroids on a singleton base space are Lie algebras. Another
extreme example is the tangent bundle $TM$ with the canonical bracket
on the space $\X(M)=\zG TM$ of vector fields.

Lie pseudoalgebras in slightly more general setting  appeared  first  in
the paper of Herz \cite{He} but one can find similar concepts under more
than a dozen  of  names  in the   literature   (e.g.   $(R,A)$-Lie
algebras,   Lie-Cartan   pairs, Lie-Rinehart algebras, differential
algebras, etc.). Lie algebroids were introduced  by Pradines \cite{Pr}.
For both notions we refer to a survey article  by Mackenzie \cite{Ma}.

From  now  on  we  assume  that  $L$  is  a  vector  bundle  over   $M$,
$V(M)=C^\infty(M)$ is the algebra of smooth functions on $M$, $\Li$ is  the
$V(M)$-module of smooth sections of $L$. Any Gerstenhaber algebra  structure
on the Grassmann algebra $\A(L)=\oplus_{n\in\Z}\A^n(L)$, where
$\A^n(L)=\zG(\bigwedge^{n+1}L)$, we will  call  a  {\em  Schouten-Nijenhuis
algebra}. As it was already indicated in \cite{KS},  Schouten-Nijenhuis
algebras are in one-one correspondence with Lie algebroids:

\begin{theorem} Any  Schouten-Nijenhuis  bracket $[\cdot,\cdot]$  on
$\A(L)$ induces a Lie  algebroid  bracket  on  $\Li=\A^0(L)$  with  the
anchor defined by $a(X)(f)=[X,f]$. Conversely, any Lie algebroid
structure  on
$\Li$ gives rise to a  Schouten-Nijenhuis  bracket  on  $\A(L)$  for  which
$\Li=\A^0(L)$ is a Lie subalgebra and $a(X)(f)=[X,f]$.
\end{theorem}
Let   $\zW(L)=\oplus_{n\in\Z}\zW^{[n]}$ be the $V(M)$-module   dual   to
$\A(L)=\oplus_{n\in\Z}\A^{[n]}(L)$,                                where
$\zW^{[n]}(L)=\zG(\bigwedge^nL^*)$  is  the
space of sections of the $n$-exterior power of the bundle $L^*$ dual  to
$L$. We can think of elements of $\A^{[n]}(L)$ as being `$n$-vector fields'
and elements of $\zW^{[n]}(L)$ as being `n-forms'.

The  Lie  algebroid  bracket  on  $\Li=\zG  L$  induces  the  well-known
generalization of the standard Cartan calculus of differential forms and
vector fields \cite{Ma,MX} (one can find an algebraic calculus for gauge
theories built on extensions of Lie algebroids in \cite{LM}).

The exterior derivative $\D:\zW^{[k](L)} \rightarrow\zW^{[k+1]}(L)$
is defined by the standard formula
\bea\nn
        \D\zm(X_1,\dots,X_{k+1}) &=& \sum_i (-1)^{i+1}
[X_i,\zm(X_1,\dots,\widehat{X}_i,\dots ,X_{k+1})] \\
&+& \sum_{i<j} (-1)^{i+j}\zm ([X_i,X_j], X_1,\dots , \widehat{X}_i,
\dots ,\widehat{X}_j, \dots , X_{k+1}),\label{D}
\eea
        where $X_i\in\Li$ and the hat over a symbol means
that this is to be omitted.   For $X \in\Li$, the
contraction $i_X \colon \zW^{[p]}(L) \rightarrow\zW^{[p-1]}(L)$
is defined in the standard way and the Lie differential operator
$ \Ll_X$
         is defined by the graded commutator
\be\label{L}
\Ll_X = i_X \circ \D - \D \circ i_X.
\ee
        The following theorem contains a list of well-known
properties of these objects.
\begin{theorem}
        Let  $\zm\in \zW^{[k]}(L),\ \zn\in \zW(L))$ and $X,Y \in\Li$.
We have
        \begin{description}
        \item{(a)} $\D \circ \D =0$,
        \item{(b)} $\D(\zm\wedge \zn) = \D(\zm) \wedge \zn +
(-1)^k\zm\wedge \D(\zn)$,
        \item{(c)} $\Ll_X(\zm\wedge \zn) = \Ll_X\zm\wedge \zn +
\zm\wedge \Ll_X \zn$,
        \item{(d)} $\Ll_X\circ \Ll_Y - \Ll_Y \circ \Ll_X =\Ll_{[X,Y]}$,
        \item{(e)} $\Ll_X \circ i_Y - i_Y\circ \Ll_X = i_{[X,Y]}$.
        \end{description}
\end{theorem}

Let us now consider Gerstenhaber-Jacobi  structures  on  the
Grassmann  algebra  $\A(L)$.  We  will  call  these   {\em   Schouten-Jacobi
algebras}.  In  view  of  the  previous  theorem,   we   will   identify
Schouten-Jacobi brackets on $\A(L)$ with {\em Jacobi algebroid} structures
on $L$. We will see that a Jacobi algebroid structure on  $L$
is determined by a  Lie  algebroid  structure  on  $L$  and  a  `1-form'
$\zF\in\zW^{[1]}(L)$   which   is   closed,   i.e.    $\D\zF=0$.    Indeed,
since $D(X)=(-1)^x\tilde D(X)$ defines a graded derivative of the  wedge
product  of  degree  -1,  $D$  is  $V(M)$-linear,  so  $D=i_\zF$  for  some
$\zF\in\zW^{[1]}(L)$. Moreover,  since  $\tilde  D(X)=[X,\1]$,  the  graded
Jacobi identity implies that $D$  is  a  derivative  also  for  the  Lie
bracket, i.e. $D([X,Y])=[D(X),Y]+[X,D(Y)]$ for $X,Y\in\A^0(L)$ which  means
exactly that $\D\zF=0$. Further, the Schouten-Jacobi bracket  restricted
to $\A^0(L)$ is a Lie algebroid bracket. Indeed, (\ref{F}) implies in
particular that $[X,fg]=[X,f]g+f[X,g]-D(X)fg$  for  $X\in\A^0(L)$,  $f\in
V(M)$, which means that $\ad_X$ induces on $V$ a first  order  differential
operator  $\hat X+D(X)I$,  where  $\hat  X$   is   a   derivation   and
$D(X)=\ad_X(\1)$ is the  part  of  order  0  (cf. Theorem 2).  Since  for
$X,Y\in\A^0$, $f\in V$, we have in view of (\ref{F})
\be
[X,fY]=([X,f]-D(\1)f)Y+f[X,Y]=\hat X(f)Y+f[X,Y],
\ee
the bracket on  $\A^0(L)$ is a  Lie  algebroid  bracket  with  the  anchor
$a(X)=\hat X$.
\par
Conversely,  having  a  Lie   algebroid   bracket   $[\cdot,\cdot]$   on
$\A^0(L)$   with an   anchor   $a$    and    a     closed     `1-form'
$\zF\in\zW^{[1]}(L)$  we
can construct a Schouten-Jacobi bracket  on  $\A(L)$  as  it  was  done  in
\cite{IM} (but here the signs are adapted to our conventions):
\be\label{IM}
[X,Y]^\zF=[X,Y]^{SN}+xX\we i_\zF Y-(-1)^xyi_\zF X\we Y,
\ee
where  $[\cdot,\cdot]^{SN}$  is  the   Schouten-Nijenhuis   bracket   on
$\A(L)$ induced by the Lie algebroid structure and $x,y$ are Lie algebra
degrees of $X,Y$. This bracket gives the  original
Lie algebroid structure on $\A^0$ and $\tilde D(X)=[X,\1]=(-1)^xi_\zF X$
for $X\in\A(L)$. Thus we get the following.
\begin{theorem}   The   formula   (\ref{IM})   describes    a    one-one
correspondence between  Schouten-Jacobi  brackets  on  $\A(L)$  and  Jacobi
algebroid  structures on   $L$,   i.e.   Lie   algebroid   brackets    on
$\A^0(L)$
defining Schouten-Nijenhuis brackets $[\cdot,\cdot]^{SN}$ with the presence
of a 1-cocycle $\zF\in\zW^{[1]}(L)$, $\D\zF=0$.
\end{theorem}
One can develop a Cartan calculus for Jacobi algebroids similarly to the
Lie algebroid  case  (cf.  \cite{IM}).  For  a  Schouten-Jacobi  bracket
associated with  a  1-cocycle  $\zF$  the  definitions  of  the  exterior
differential  $\D^\zF$   and   Lie   differential   $\Ll^\zF=\D^\zF\circ
i+i\circ\D^\zF$  are formally the same as (\ref{D}) and
(\ref{L}),  respectively.   Since,   for $X\in\A^0(L)$, $f\in V(M)$,   we
have $[X,f]=[X,f]^{SN}+(i_\zF X)f$, one obtains $\D^\zF\zm=\D\zm+
\zF\we\zm$. Here  $[\cdot,\cdot]^{SN}$  and  $\D$
are, respectively,  the  Schouten-Nijenhuis  bracket  and  the   exterior
derivative associated with the Lie algebroid of the given Jacobi algebroid.
For the exterior differential and Lie differential associated with a Jacobi
algebroid we have the following.
\begin{theorem}
        Let  $\zm\in \zW^{[k]}(L),\ \zn\in \zW(L))$ and $X,Y \in\Li$.
We have
        \begin{description}
        \item{(a)} $\D^\zF \circ \D^\zF =0$,
        \item{(b)} $\D^\zF(\zm\wedge \zn) = \D^\zF(\zm) \wedge \zn +
(-1)^k\zm\wedge \D^\zF(\zn)-\zF\we\zm\we\zn$,
        \item{(c)} $\Ll^\zF_X(\zm\wedge \zn) = \Ll^\zF_X\zm\wedge \zn +
\zm\wedge \Ll^\zF_X \zn -(i_\zF X)\zm\we\zn$,
       \item{(d)}  $\Ll^\zF_X  \circ  i_Y   -   i_Y\circ   \Ll^\zF_X   =
i_{[X,Y]}$,
        \item{(e)} $\Ll^\zF_X\circ \Ll^\zF_Y - \Ll^\zF_Y \circ \Ll^\zF_X
=\Ll^\zF_{[X,Y]}$.
        \end{description}
\end{theorem}
\begin{pf} The proof of (a), (b), (c) can be found  in  \cite{IM}.  The
property (d) follows easily from definitions,  and  (e)  follows  easily
from (d).
\end{pf}

\noindent
{\bf Example.} Consider the Jacobi algebroid structure  associated  with
the Nijenhuis-Richardson bracket on $\A Diff_1(M)$.
According  to Theorem 2,  $\A  Diff_1^0(M)$  can  be   identified   with
$\X(M)\oplus C^\infty(M)$, i.e. with sections of the direct  sum  bundle
$TM\oplus_M\R$. It is  easy  to  see that the Lie algebroid bracket on
this bundle reads (cf. \cite{IM})
\be
[(X,f),(Y,g)]=([X,Y],X(g)-Y(f)),
\ee
where the right-hand side bracket is  the  standard  bracket  of  vector
fields. The 1-cocycle $\zF$, written as $i_\1$ in Theorem 2, is given by
$\zF(X,f)=f$. The Schouten-Jacobi bracket (i.e. the Nijenhuis-Richardson
bracket in this case) reads
\bea\nn
&&[A_1+I\we A_2,B_1+I\we B_2]^{RN}=[A_1,B_1]^{SN}+(-1)^aI\we[A_1,B_2]^{SN}
+I\we[A_2,B_1]^{SN}\\
&&\qquad+aA_1\we B_2-(-1)^abA_2\we B_1+(a-b)I\we A_2\we B_2.\label{dif}
\eea
Hence, the bracket  $\{\cdot,\cdot\}$  on  $C^\infty(M)$  defined  by  a
bilinear differential operator $\zL+I\we\zG\in\A  Diff_1^1(M)$  is  a  Lie
bracket (Jacobi bracket on $C^\infty(M)$) if and only if
\be
[\zL+I\we\zG,\zL+I\we\zG]^{RN}=[\zL,\zL]^{SN}+2I\we
[\zG,\zL]^{SN}+2\zL\we\zG=0.
\ee
We recognize the conditions
\be\label{js}
[\zG,\zL]^{SN}=0, \quad [\zL,\zL]^{SN}=-2\zL\we\zG,
\ee
defining a Jacobi structure on $M$ (\cite{Li2}). The  difference  in the
sign when comparing with \cite{Li2} comes from  our  convention  for  the
Schouten bracket.

It is obvious that the formula  (\ref{dif})  defines  a  Schouten-Jacobi
bracket  for  any  extension  $L\oplus_M\R$  of  a  Lie  algebroid  $L$
associated with the anchor map:
\be\label{an}
[(X,f),(Y,g)]=([X,Y],a(X)(g)-a(Y)(f)).
\ee
The 1-cocycle is in this case $\zF((X,f))=f$.

\medskip
Suppose  that  we  have  a  Schouten-Jacobi   bracket   as   above   and
$Z\in\A^0(L)$. We call $X\in\A^x(L\oplus\R)$  a  {\it  $Z$-homogeneous}
element   if   $[Z,X]^\zF=-xX$.   The   graded   subspace   spanned   by
$Z$-homogeneous  elements  is  clearly   a   Lie   subalgebra   of   the
Schouten-Jacobi bracket. We can represent  the  Schouten-Jacobi  bracket
for $Z$-homogeneous elements of  $\A(L\oplus\R)$  in  the
Schouten-Nijenhuis bracket of $\A(L)$.
\begin{theorem} Let  $H_Z$  be  the  mapping  which  associated with  any
$Z$-homogeneous element $A=A_1+I\we A_2\in\A^a(L\oplus\R)$  the  element
$H_Z(A)=A_1+Z\we  A_2$.  Then   $H_Z$   is   a   homomorphism   of   the
Schouten-Jacobi bracket (\ref{dif}) on $Z$-homogeneous elements into the
Schouten-Nijenhuis bracket on $\A(L)$:
\bea\nn
[H_Z(A),H_Z(B)]^{SN}&=&
[A_1+Z\we A_2,B_1+Z\we B_2]^{SN}=\\  \label{dif1}
&&[A_1,B_1]^{SN}+aA_1\we B_2-(-1)^abA_2\we B_1\\ \nn
&&+(-1)^aZ\we[A_1,B_2]^{SN}+Z\we[A_2,B_1]^{SN}+(a-b)Z\we A_2\we
B_2\\ \nn
&=&H_Z([A,B]^\zF).
\eea
\end{theorem}

\section{Lifts of Schouten-Jacobi brackets}
Since in (\ref{IM}) we can put $h\zF$ instead of $\zF$, where  $h$  is  a
parameter, the Schouten-Jacobi bracket can be viewed as a deformation of
the   Schouten-Nijenhuis   bracket.
 \begin{theorem}
The    Schouten-Jacobi    brackets $[\cdot,\cdot]^{SN}$ and
$[\cdot,\cdot]^\zF_0$ defined by
\be
[X,Y]_0^\zF=xX\we i_\zF Y-(-1)^xyi_\zF X\we Y
\ee
are compatible,i.e. $[X,Y]^{h\zF}=[X,Y]^{SN}+h[X,Y]^\zF_0$ is  a  graded
Lie bracket for all $h\in\R$.
\end{theorem}
An easy way to see that the Schouten-Jacobi bracket (\ref{IM}) is really a
graded Lie bracket is the following.
\par
Consider the product $\tilde L$ of the Lie algebroid $L$ and $T\R$, i.e.
we view $\tilde L=L\ti T\R$ as a vector bundle  over  $M\ti\R$  with  the
obvious product Lie bracket and the  anchor  $a\ti{\rm  id}:L\ti  T\R\ra
TM\ti T\R$. For a fixed 1-cocycle $\zF\in\zW^{[1]}(L)$ we define  a  Lie
algebroid injective homomorphism $U_\zF:L\ra\tilde L$ by
\be
U_\zF(X)=X+i_\zF X\pa_t.
\ee
Here sections of $L$ and functions on $M$  on  the  right-hand  side  are
understood as sections of  $L\ti  T\R$  and  functions  on  $M\ti\R$  in
obvious way and $\pa_t$ is the basic vector field on $\R$. Since $\zF$  is
a 1-cocycle,
\be
[U_\zF(X),U_\zF(Y)]=[X,Y]+([X,i_\zF Y]+[i_\zF X,Y])\pa_t=[X,Y]+i_\zF
[X,Y]\pa_t=
U_\zF([X,Y])
\ee
for $X,Y\in\Li$ and $U_\zF$ is really a homomorphism. This homomorphism can be
extended to  a homomorphism of the whole Gerstenhaber algebra by
\be\label{U}
U_\zF(X)=X+\pa_t\we i_\zF X,
\ee
since $U_\zF$ respects the wedge product. Thus
\be
[U_\zF(X),U_\zF(Y)]^{SN}=U_\zF([X,Y]^{SN})
\ee
for all $X,Y\in\A(L)$.
Now, we can define  a  new  graded
linear    map    $\tilde    U_\zF:\A(L)\ra\A(\tilde    L)$    by     $\tilde
U_\zF(X)=e^{-xt}U_\zF(X)$. This mapping respects  grading  but  not  the
wedge
product. However, the image of $\tilde U_\zF$ is a  Lie  subalgebra  of  the
Schouten-Nijenhuis bracket on $\A(\tilde L)$. It  is  easy  to  see  the
following.
\begin{theorem} For $X,Y\in\A(L)$ we have
\be
[\tilde U_\zF(X),\tilde U_\zF(Y)]^{SN}=\tilde U_\zF([X,Y]^{SN}+xX\we i_\zF
Y-(-1)^xyi_\zF X\we Y).
\ee
Thus   $\tilde   U_\zF$   is   an   embedding   of   the     Schouten-Jacobi
bracket $[\cdot,\cdot]^\zF$  on $\A(L)$ (induced by   the
Schouten-Nijenhuis bracket and the 1-cocycle
$\zF$)  into  the  Schouten-Nijenhuis  bracket  on  $\A(\tilde  L)$,  so
(\ref{IM}) is a graded Lie bracket.
\end{theorem}
Note that on a similar idea is based the Poissonization of a
Jacobi structure in  \cite{GL}  and the construction a Lie algebroid
from a Jacobi structure in \cite{Va1}.
The bundle projection of $L\ti T\R$ over  $M\ti\R$  onto  $L\ti\R$  over
$M\ti\R$  defines  a  Lie  algebroid bracket  on   $L\ti\R$,   i.e.
on `time-dependent sections of $L$' as described  in  \cite{IM}.
Composing $\tilde U_\zF$ with this projection we  get  a  representation  of
the  Schouten-Jacobi  bracket  $[\cdot,\cdot]^\zF$  on  $\A(L)$  in  the
Schouten-Nijenhuis bracket of the Lie algebroid $L\ti\R$ (cf. \cite{IM},
section 4.2). The advantage of this construction is that  the  dimension
of the fibres remains the same. On the other hand, in  our  construction
the Lie algebroid is fixed and only the embedding depends on $Phi$.
\par
There is another  approach  to  Lie  algebroids.  As  it  was  shown  in
\cite{GU1, GU2}, a Lie algebroid   structure   (or   the   corresponding
Schouten-Nijenhuis  bracket)  is  determined  by  the   algebroid   lift
$X\mapsto X^c$ which associates with  $X\in\A(L)$  a  multivector  field
$X^c\in A(M)$. Recall, that sections $\zm$ of the dual bundle $L^*$  may
be identified with linear (along fibres)  functions  $\zi_\zm$  on  $L$:
$\zi_\zm(X_p)=<\zm(p),X_p>$.  By  {\it  homogeneous  elements}  of   the
Schouten algebra of multivector fields on a vector bundle we  understand
elements which are homogeneous with  respect  to  the  Liouville  vector
field     $\nabla$.     This     means      that each contraction
with differentials of linear functions
$<\zL,\D\zi_{\zm_{0}}\we\cdots\we\D\zi_{\zm_{\zl}}>$ is again a linear
function associated with  an  element  $[\zm_0,\dots,\zm_\zl]_\zL$.  The
multilinear operation $[\zm_0,\dots,\zm_\zl]_\zL$ on sections  of  $L^*$
we call the {\it bracket induced by $\zL$}. Note  that  these  brackets
have a property similar to the Lie algebroid brackets:
\be
[\zm_0,\dots,\zm_{\zl-1},f\zm_\zl]=f[\zm_0,\dots,\zm_{\zl-1},\zm_\zl]+
\zL_{\zm_0,\dots,\zm_{\zl-1}}(f)\zm_\zl,
\ee
where
\be
\zL_{\zm_0,\dots,\zm_{\zl-1}}(f)=[\zm_0,\dots,\zm_{\zl-1},\D f]_\zL
\ee
defines the {\it  Hamiltonian  vector field} of $\zL$ associated with
$(\zm_0,\dots,\zm_{\zl-1})$. Thus, there  is  a  one-one  correspondence
between linear multivector fields and such brackets.

\begin{theorem} (\cite{GU1}) For a given Lie  algebroid  structure  on  a
vector bundle $L$ over $M$ there is a unique {\it complete lift} of
elements $X$ of the Gerstenhaber algebra $\A(L)$ to homogeneous  elements
$X^c$  of the Schouten algebra $A(L)$ of multivector fields on $L$, such
that
\begin{description}
\item{(a)} $f^c=\zi_{\D f}$ for $f\in C^\infty(M)$;
\item{(b)}  $X^c(\zi_\zm)=\zi_{\Ll_X\zm}$  for  $X\in\zG  L,\  \zm\in\zG
L^*$;
\item{(c)} $(X\we Y)^c=X^c\we Y^v+X^v\we Y^c$, where $X\mapsto  X^v$  is
the standard vertical lift of multisections of $L$ to multivector  fields
on $L$.
\end{description}
Moreover,  this  complete   lift   is   a   homomorphism   of   the
Schouten-Nijenhuis brackets:
\be\label{l1}
[X,Y]^c=[X^c,Y^c]
\ee
and
\be\label{l2}
[X^c,Y^v]=[X,Y]^v.
\ee
\end{theorem}
{\bf Remark.} For the canonical Lie algebroid $L=TM$, the above complete
lift reduces to the better-known {\it tangent lift} of multivector  fields
on $M$ to multivector fields on $TM$ (cf.  \cite{IY,GU}).  The  complete
Lie algebroid lift of just sections of $L$, i.e. the  formula  (b),  was
already indicated in \cite{MX1}.

\medskip\noindent
Our aim is to find an analog of the  Lie  algebroid  complete  lift  for
Jacobi algebroids which will represent the  Schouten-Jacobi bracket  on
$\A(L)$    in    the    Nijenhuis-Richardson    bracket     of     first
order multidifferential operators on $L$. Let $[\cdot,\cdot]^\zF$ be the
Schouten-Jacobi bracket on  $\A(L)$  associated  with  a  Lie  algebroid
structure on $L$ and a 1-cocycle $\zF$.

\medskip\noindent
{\bf Definition.} The  {\it  Jacobi  lift}  of  an
element $X\in\A^x(L)$ is the  element  $\hat X_\zF\in\A  Diff_1(L)$,
i.e.  a multidifferential operator of first order on $L$, defined by
\be\label{jl}
\hat X_\zF=X^c-x\zi_\zF X^v+I\we(i_\zF X)^v,
\ee
where $X^c$ is the complete Lie algebroid lift and $X^v$ is the  vertical
lift.
\begin{theorem} The Jacobi lift has the following properties:
\begin{description}
\item{(a)} $\hat f_\zF=\zi_{\D^\zF f}$ for $f\in C^\infty(M)$;
\item{(b)}   $\hat X_\zF(\zi_\zm)=\zi_{\Ll^\zF_X\zm}$ and
$X_\zF(\1)=\zF(X)\circ\tau$
for
$X\in\zG L,\ \zm\in\zG L^*$ and $\tau:L\rightarrow M$ being  the  bundle
projection;
\item{(c)} $(X\we Y)^{\hat{}}_\zF=\hat X_\zF\we Y^v+X^v\we
\hat Y_\zF-\zi_\zF(X^v\we Y^v)$;
\item{(d)} $[\hat X_\zF,\hat Y_\zF]^{NR}=([X,Y]^\zF)^{\hat{}}_\zF$.
\end{description}
\end{theorem}
\begin{pf}  The  proof  consists  of  standard  calculations  using  the
properties of the Schouten-Nijenhuis and  Schouten-Jacobi  brackets  and
the properties of
the complete lift. One should  also  remember  that  $i_\zF[X,Y]^\zF
=[i_\zF X,Y]^\zF+(-1)^x[X,i_\zF Y]^\zF$ and  use  the  identities   (cf.
\cite{GU1})
$[\zi_\zF,X^v]^{SN}=-(i_\zF X)^v$ and $[\zi_\zF,X^c]^{SN}=-(i_\zF  X)^c$
(the last one depends on $\D\zF=0$).
\end{pf}

\begin{corollary} If $X=\zL+I\we\zG$ is a Jacobi structure on $M$,  then
the  Jacobi  lift  $\hat X_\zF$  is  a  homogeneous   Jacobi structure   on
$TM\oplus\R$. Moreover,
\be
\hat X_\zF=\zL^c+\pa_t\we\zG^c-t(\zL^v+\pa_t\we\zG^v)+I\we\zG^v,
\ee
where $\zL^c$, $\zL^v$, etc., are the complete  and  vertical  lifts  to
$TM$ and $t$ is the canonical linear coordinate in $\R$.
\end{corollary}
\begin{pf} In our case the Lie algebroid is the  extension  $TM\oplus\R$
relative to the anchor map and the complete and vertical lifts of  $\zL$
and $\zG$ with respect to this Lie algebroid structure coincide with the
standard  tangent  complete  and  vertical  lifts.  Moreover,   $I^c=0$,
$I^v=\pa_t$ and $i_\zF X=\zG$.
\end{pf}

\noindent
{\bf  Remark.}  We  get  the  same  homogeneous  Jacobi   structure   as
\cite{IM1}, example 5.

\medskip\noindent
According to Theorem 7,  there  is  a  homomorphism  $H_\nabla:\A
Diff_1(L)\ra\A(L)$ of  the  Nijenhuis-Richardson  bracket  on  homogeneous
first order multi-differential operators into the Schouten-Nijenhuis
bracket  of  multivector   fields on  $L$   given   by   $H_\nabla(X_1+I\we
X_2)=X_1+\nabla\we    X_2$, where $\nabla$ is the Liouville vector field
on the vector bundle $L$.

\medskip\noindent
{\bf  Definition.}  The      {\it      Poisson      lift}  $X^c_\zF$  is
defined by
\be\label{pl} X^c_\zF=H_\nabla(\hat X_\zF)=
X^c-x\zi_\zF X^v+\nabla\we(i_\zF X)^v.
\ee
\begin{theorem} The Poisson lift has the following properties:
\begin{description}
\item{(a)} $f^c_\zF=\zi_{\D^\zF f}$ for $f\in C^\infty(M)$;
\item{(b)}   $X^c_\zF(\zi_\zm)=\zi_{\Ll^\zF_X\zm}$;
\item{(c)} $(X\we Y)^c_\zF=X^c_\zF\we Y^v+X^v\we Y^c_\zF-\zi_\zF(X^v\we
Y^v)$;
\item{(d)} $[X^c_\zF,Y^c_\zF]^{SN}=([X,Y]^\zF)^c_\zF$.
\end{description}
\end{theorem}
\section{Lie algebroids associated with local Lie algebras}

It  is  known  (\cite{Fu,Ko,GU}), that a Poisson  structure  $\zL$ on
$M$ defines  not  only
the  Poisson bracket $\{\cdot,\cdot\}_\zL$ of functions, but also a Lie
bracket $[\cdot,\cdot]_\zL$ on 1-forms, given by
\be
[\zm  ,\zn]_\zL = \Ll_{\zL_\zm}\zn   -\Ll_{\zL_\zn}\zm   -\D
<\zL,\zm\we\zn>,\label{P8}
\ee
where $\zL_\zm=i_\zm\zL$ and $<\cdot ,\cdot >$ is  the  pairing  between
forms  and multivector fields.
In particular, $[\D f,\D g]_\zL =\D\{ f,g\}_\zL$ and $\zL^\#$  is  a  Lie
bracket homomorphism:
\be
[\zL_\zm ,\zL_\zn ]_\zL=\zL_{[\zm ,\zn]_\zL}.\label{P9}
\ee
This bracket on 1-forms is a Lie algebroid bracket and  it  induces  the
corresponding Schouten-Nijenhuis bracket on $\zW (M)$.
It was observed by Koszul \cite{Ko} (see also   \cite{KSM})  that  this
bracket has a {\it generating operator } $\pa_\zL$:
\be
[\zm ,\zn]_\zL  =(-1)^m  (\partial_\zL(\zm\we\zn  )-\partial_\zL\zm\we
\zn -(-1)^m \zm\we\partial_\zL\zn),\label{P10}
\ee
where  $\partial_\zL=\I  (P)\circ\D  -\D\circ\I  (P)$  and  $m$  is  the
standard degree of the form $\zm$. Note  that $\pa_\zL$ is a homology
operator, since $\pa_\zL^2=0$.Moreover,
\be
\D [\zm ,\zn]_\zL =[\D\zm ,\zn]_\zL +(-1)^{m -1}[\zm ,\D\zn]_\zL.
\label{P12}
\ee
The whole structure, i.e. the Gerstenhaber algebra with  the  generating
operator $\pa$ for the Lie bracket satisfying $\pa^2=0$  is  called
 a  {\it Batalin-Vilkovisky algebra}  (cf.  \cite{KS}).  With the  presence
of  the
derivation $\D$ of both: the associative and Lie algebra structures it is
a {\it differential Batalin-Vilkovisky algebra}.

It is well-known \cite{GU} that the Lie algebroid bracket (\ref{P8})  is
induced by the complete lift $\zL^c$. Now, it should be no surprise  that
the Lie algebroid bracket (\ref{jb}) we started with is induced  by  our
Poisson lift of the corresponding Jacobi structure.

\begin{theorem} If $X=\zL+I\we\zG$ is a Jacobi structure on $M$,  then
the Poisson  lift  $X^c_\zF$  is  a  homogeneous  Poisson   structure   on
$TM\oplus\R$. Moreover,
\be
X^c_\zF=\zL^c+\pa_t\we\zG^c-t\zL^v+\nabla\we\zG^v,
\ee
where $\zL^c$, $\zL^v$, etc., are the complete  and  vertical  lifts  to
$TM$, respectively, $\nabla$ is the Liouville vector field on $TM$ and
$t$ is the canonical linear coordinate in
$\R$.   This homogeneous Poisson structure determines a Lie  algebroid
bracket  on the dual bundle $T^*M\oplus\R$, given for
$\zm,\zn\in\zG(T^*M\oplus\R)$,
$\zm=(\za,f)$, $\zn=(\zb,g)$, by
\be\label{br}
[\zm,\zn]=\Ll_{X_\zm}^\zF\zn-\Ll_{X_\zn}^\zF\zm-\D^\zF<X,\zm\we\zn>,
\ee
where $X_\zm=i_\zm X=(i_\za\zL+f\zG,-i_\za\zG)$, $X_\zn=i_\zn X=
(i_\zb\zL+g\zG,-i_\zb\zG)$ are the  corresponding  first  order  operators
viewed as sections of $TM\oplus\R$ and $\zF$ is the 1-cocycle defined by
the projection on $\R$.
This bracket coincides with (\ref{jb}).
\end{theorem}
\begin{pf} The first part follows immediately from the  general  result.
To   show (\ref{br}) consider $X=X_1\we X_2$ and
$\zm=\zm_1\we\zm_2$, where $X_i\in\A^0$, $\zm_i\in\zW^{[1]}$.  Then  the
bracket $\{\zi_{\zm_1},\zi_{\zm_2}\}_{X^c_\zF}$ defined by $X^c_\zF$ is
given by
\be
<(X_1)^c_\zF\we (X_2)^v+(X_1)^v\we(X_2)^c_\zF-\zi_\zF(X_1^v\we
X_2^v),\D\zi_{\zm_1}\we\D\zi_{\zm_2}>.
\ee
This is the linear function
\be
\sum_{i,j}(-1)^{i+j}{(X_i)}^c_\zF(\zi_{\zm_j})<X_{\zs(i)}^v,
\zi_{\zm_{\zs(j)}}>-\zi_\zF<X^v,\D\zi_{\zm_1}\we\D\zi_{\zm_2}>,
\ee
where $\zs$ is the transposition of $(1,2)$.
Taking  into  account  that  $Y^c_\zF(\zi_\zn)=\zi_{\Ll^\zF_Y\zn}$  and
$Y^v(\zi_\zn)=\zi_{<Y,\zn>}$ we get that this linear function corresponds
to
\be
\sum_{i,j}(-1)^{i+j}<X_i,\zm_j>\Ll^\zF_{X_{\zs(i)}}\zm_{\zs(j)}-<X,\zm>\zF.
\ee
In view of properties $f\Ll^\zF_Y\zn=\Ll^\zF_{fY}\zn-<Y,\zn>\D^\zF f$ and
$\D^\zF(fg)=f\D^\zF g+g\D^\zF f-fg\zF$, we get (\ref{br}). To show  that
this is exactly (\ref{jb}) it suffices to lead  easy  calculations  with
$\Ll^\zF_{Y+fI}\zn=\Ll_Y\zn+f\zn$.
\end{pf}

\noindent
Of course, instead of the Jacobi structure $X=\zL+I\we\zG$  we  can  start
from an arbitrary Schouten-Jacobi bracket $[\cdot,\cdot]^\zF$ on $\A(L)$
and an element $X\in\A^1(L)$ with $[X,X]^\zF=0$ (we will call such $X$ a
{\it Jacobi element}) and use the formula (\ref{br})  to  define  a  Lie
algebroid structure on $L^*$ (see \cite{IM}).

\noindent
{\bf Example.} Let $L$ be a Lie algebroid and let $(\zL,\zG)$ be a  {\it
Lie    algebroid    Jacobi    structure},     i.e.      $\zL\in\A^1(L)$,
$\zG\in\A^0(L)$, satisfy (\ref{js}). Then, the  element  $X=\zL+I\we\zG$
regarded  as  a bisection of the Lie algebroid extension  $L\oplus_M\R$
relative  to  the anchor map (cf. (\ref{an})) is a Jacobi element  of
the  Schouten-Jacobi bracket  associated  with  the   1-cocycle  given
by $\zF((Y,f))=f$ (cf. example 1).  Thus,  the  Poisson  lift  $X^c_\zF$
defines a Lie algebroid bracket on $L^*\oplus_M\R$ which is formally  the
same as (\ref{jb}).

The  formula  (\ref{br})  for  the  bracket  of  `1-forms' generated by
$X^c_\zF$ can be generalized for arbitrary $X\in\A^x(L)$ as follows.
\begin{theorem} The bracket $[\zm_0,\dots,\zm_x]_{X^c_\zF}$
for $\zm_i\in\zG  L^*$ induced by $X^c_\zF$ by
\be
\zi_{[\zm_0,\dots,\zm_x]_{X^c_\zF}}=
<X^c_\zF,\D^\zF\zi_{\zm_0}\we\cdots\we\D^\zF\zi_{\zm_x}>
\ee
reads
\be\label{o}
[\zm_0,\dots,\zm_x]_{X^c_\zF}=\sum_{k=0}^x(-1)^{x+k}\Ll_{X_k}\zm_k-x\D^\zF
<X,\zm_0\we\cdots\we\zm_x>,
\ee
where
\be
X_k=i_{\zm_0\we\cdots\we\check\zm_k\we\cdots\we\zm_x}X.
\ee
In particular,
\be\label{lo}
[\D^\zF f_0,\dots,\D^\zF  f_x]_{X^c_\zF}= \D^\zF<X,\D^\zF   f_0\we \cdots
\we\D^\zF f_x>.
\ee
\end{theorem}
{\bf Remark.} In the case when the  differentials  $\D^\zF  f$  generate
$L^*$ almost everywhere the formula (\ref{lo}) defines the bracket, thus
the  homogeneous tensor  $X^c_\zF$, uniquely.

\medskip\noindent
Suppose now  that  $L$  is  a one-dimensional  vector  bundle  over  $M$.
Like in the case of the trivial bundle, the graded space $\A Diff_1(\Li)$
of first order multidifferential operators on $L$ is a Lie subalgebra of
the Richardson-Nijenhuis bracket on $\A(\Li)$ -- the space of multilinear
maps of the vector space $V=\Li$ of sections of $L$. The difference with
the case $V=C^\infty(M)$ is that we do not have  a  natural  associative
algebra structure on $V$. However locally, fixing a basic  section  over
an  open  subset   $N$  of  $M$,  we  have  an   isomorphism   of    the
corresponding  graded
Lie algebras $\A Diff_1(\Li_{\mid N})$ and $\A Diff_1(N)$. Thus, locally, we
have the Poisson lift $A^c_\zF$ for any $A\in\A Diff_1(\Li)$. The problem is
that  what  corresponds  to  $\1$,  and  hence  what  is  $\zF$  on  $\A
Diff_1(\Li)$,
depends on the choice of the local section. However, using the preceding
remark, we can conclude that there is a uniquely  defined  complete  lift
$A^c_{loc}$. This time it is not a multivector field  on  $TM\oplus_M\R$
but on the dual $J_1^*(L)$ to the first jet bundle $J_1(L)$. Of course,
in  the trivial case, $J_1(L)$ is canonically isomorphic with
$T^*M\oplus\R$, so $J_1^*(L)$ is isomorphic  with  $TM\oplus\R$  but  for
non-trivial bundles it is not the case. We just define the complete lift
$A^c_{loc}$ as the unique homogeneous $(a+1)$-vector field  on  $J^*_1(L)$
such that
\be\label{loc}
<A^c_{loc},\D \zi_{j_1(f_0)}\we\dots\we\D \zi_{j_1(f_a)}>=
j_1(A(f_0,\dots,f_a))
\ee
for all sections $f_0,\dots,f_a$ of $L$.  Here  $j_1$  means  the  first
jet prolongation of a given section. Since the first  jet  prolongations
of  sections  of  $L$  generate  $J_1(L)$ over an open-dense  subset,   the
multivector field $A^c_{loc}$ is uniquely defined. Moreover, it is  easy
to see that for any  local  trivialization  $A^c_{loc}$  coincides  with
$A^c_\zF$. Since all our  brackets  are  local  over  $M$,  we  get  the
following.

\begin{theorem} The complete lift of first order differential  operators
on a one-dimensional bundle defined by (\ref{loc}) is  a  homomorphism  of
the   Nijenhuis-Richardson   bracket   on   $\A Diff_1(\Li)$    into    the
Schouten-Nijenhuis  bracket  of  homogeneous   multivector   fields   on
$J_1^*(L)$.
\end{theorem}
\begin{corollary} If $X\in\A^1Diff_1(\Li)$ represents a local  Lie  algebra
bracket on $L$, then $X^c_{loc}$ induces a Lie  algebroid  structure  on
the first jet bundle $J_1(L)$.
\end{corollary}

\noindent
{\bf Remark.} There is no clear analog of the Jacobi lift for
$\A Diff_1(\Li)$,  since locally it  depends  stronger  on  the  1-cocycle
$\zF$ associated with the  trivialization.  This  suggests  that  the
lift  to multivector fields is primary with respect to the Jacobi lift.
Note  also  that   $\A Diff_1(\Li)$   is   a   $C^\infty(M)$-module   and
the Nijenhuis-Richardson bracket on $Diff_1(L)$ is a  Lie  algebroid
bracket like in the trivial case. This time however we have no splitting
$Diff_1(\Li)=Der(\Li)\oplus C^\infty(M)$ (it makes no sense) but we have
an exact sequence
\be
0\ra C^\infty(M)\ra Diff_1(\Li)\ra \X(M)\ra 0
\ee
which gives the anchor map for this Lie algebroid and splits  only  when
the bundle $L$ is trivializable.

\section{Jacobi bialgebroids}

Recall that  a  {\em  Lie  bialgebroid}  \cite{KS,MX}  is  a  dual  pair
$(L,L^*)$ of vector bundles equipped with Lie algebroid structures  such
that the differential $\D_*$ induced from the Lie algebroid structure on
$L^*$ as defined by (\ref{D}) is a derivation of the  Schouten-Nijenhuis
bracket induced by the Lie algebroid structure on $L$:
\be\label{bi}
\D_*[X,Y]=[\D_*X,Y]+(-1)^x[X,\D_*Y]\quad {\rm for all}\quad X,Y\in\A(L).
\ee
For Jacobi algebroids we will keep formally the same definition.

\medskip\noindent
{\bf Definition.}
A {\em Jacobi  bialgebroid}  is  a  dual  pair $(L,L^*)$ of vector
bundles equipped with Jacobi algebroid structures  such that the
differential $\D_*$ induced from the Jacobi algebroid structure on
$L^*$ is a derivation of the  Schouten-Jacobi bracket induced by the
Jacobi algebroid structure on $L$.

\medskip\noindent
{\bf  Remark.}  Note  that  the  above  definition  coincides  with   the
definition of {\it generalized Lie bialgebroid} in \cite{IM}.
Indeed,  the  condition  (4.1)  in  \cite{IM}  is  just  (\ref{bi})  for
$X,Y\in\A^0(L)$, and the condition (4.2) is just (\ref{bi}) for $Y=\1$,  so
every Jacobi bialgebroid is a generalized Lie algebroid in the  sense  of
\cite{IM}. To show the converse we will use the following lemma.
\begin{lemma} If $\D_*$ is the differential on  $\A$  associated  with  a
Jacobi algebroid structure on $L^*$, then
\bea\label{bi1}
&\D_*[X,Y\we Z]-[\D_*X,Y\we Z]-(-1)^x[X,\D_*(Y\we Z)]\\
&=(\D_*[X,Y]-[\D_*X,Y]-(-1)^x[X,\D_*Y])\we Z\nn\\
&+(-1)^{x(y+1)}Y\we(\D_*[X,Z]-[\D_*X,Z]-(-1)^x[X,\D_*Z])\nn\\
&-(\D_*\tilde D(X)-\tilde D(\D_*X)-(-1)^x[X,X_0])\we Y\we Z,\nn
\eea
where $X_0=\D_*\1$.
\end{lemma}
\begin{pf} The proof consists of standard calculations  using  (\ref{F})
and the following property of the exterior differential:
\be
\D_*(Y\we Z)=(\D_*Y)\we Z+(-1)^{y+1}Y\we(\D_*Z)-X_0\we Y\we Z.
\ee
\end{pf}

\begin{theorem} If (\ref{bi}) is satisfied for all $X,Y\in\A^0(L)$ and  all
$X\in V(M)\oplus\A^0(L), Y=\1$, then it is satisfied in general.
\end{theorem}
\begin{pf} First, note that
\be
\D_*\tilde D(X)-\tilde D(\D_*X)-(-1)^x[X,X_0]=
\D_*[X,\1]-[\D_*X,\1]-(-1)^x[X,\D_*\1]=0
\ee
for $X\in V(M)\oplus\A^0(L)$. Hence,  for  $X,Y\in\A^0(L)$  and
$f\in V(M)$ we get from (\ref{bi1}) and (\ref{bi}) for elements of
$\A^0(L)$:
\be
(\D_*[X,f]-[\D_*X,f]-[X,\D_*f])\we Y=0,
\ee
so (\ref{bi}) is satisfied also for $X\in\A^0(L)$ and  $Y\in  V(M)$.
Using, in turn, this fact when applying  $X=f,Y=g\in V(M)$,
$Z\in\A^0(L)$, to (\ref{bi1}), we get
\be
(\D_*[f,g]-[\D_*f,g]+[f,\D_*g])Z=0,
\ee
so (\ref{bi}) is satisfied for all $X,Y\in V(M)\oplus\A^0(L)$.
Now,  we  can prove (\ref{bi}) by induction with respect to the sum
$x+y$ of  degrees of $X$ and $Y$. If $x+y\ge 0$ and, say, $y>0$ (the case
$x=y=0$  is covered by assumption), then we can write $Y$ as a linear
combination of wedge
products $A\we B$ with  $a,b<y$  and  (\ref{bi})  follows  for  $X,Y$  by
induction in view of (\ref{bi1}).
\end{pf}

\medskip\noindent
{\bf Example.} In \cite{IM}, Theorem 5.1, it is shown  that  any  Jacobi
element $X\in\A^1(L)$, $[X,X]^\zF=0$.  for  a  Jacobi  algebroid
structure in  $\A(L)$  associated  with  a  1-cocycle  $\zF\in\A^0(L^*)$
gives rise to a Jacobi bialgebroid for which the Lie algebroid structure
on $\A(L^*)$ is given by the formula  (\ref{br})  and  the  corresponding
1-cocycle is $-X_\zF$. This is of course a Jacobi analog of a triangular
Lie bialgebroid in the sense of Mackenzie and Xu \cite{MX}.

\section{Conclusions}
We have shown that, by analogy with Lie algebroids regarded as  {\it  odd
Poisson brackets}, one can define Jacobi algebroids regarded as {\it  odd
Jacobi brackets} on  the  Grassmann  algebra  $\A(L)$  associated  with
a vector bundle $L$. Jacobi algebroids in this  sense  turn  out  to  be
objects already  studied  by  Iglesias  and  Marrero  \cite{IM}.  It  is
possible to develop a Cartan calculus for  Jacobi  algebroids.  We  have
constructed lifts of tensor fields which transport the Schouten-Jacobi
bracket on $\A(L)$ into the Schouten bracket of  multivector  fields  on
the total  space  $L$.  This  leads  to  a  natural  construction  which
associates a Lie algebroid with every local Lie algebra of Kirillov.  We
have  shown  also  that  a  notion  of  a  Jacobi  bialgebroid  can   be
consistently introduced with the full analogy to the classical case.

\medskip
Since every Lie algebroid can be viewed  as  a  particular  case  of  a
Batalin-Vilkovisky algebra \cite{Xu}, it is natural to look for   a
similar correspondence in the case of  Jacobi  algebroids. First
steps in this direction have been done in \cite{ILMP}. There  is  a
natural  way  of  defining  generating  operators  for   Schouten-Jacobi
brackets  and  of  defining  the  corresponding  homology.  We  postpone
detailed studies of these questions to a separate paper.


\end{document}